\documentclass[10pt]{amsart}
\usepackage{graphicx,color}
\newtheorem{theorem}{Theorem}[section]
\newtheorem{lemma}[theorem]{Lemma}           
\newtheorem{cor}[theorem]{Corollary}

\theoremstyle{definition}
\newtheorem{definition}[theorem]{Definition}

\theoremstyle{remark}

\numberwithin{equation}{section}

\thanks{H. Morioka was supported by the Grant-in-aid for young scientists (B) No. 16K17630, JSPS. E. Segawa was supported by the Grant-in-aid of Scientific Research (C) No. 19K036116, JSPS, and Research Origin for Dressed Photon.}

\title[Detection of edge defects by embedded eigenvalues of QW]
{Detection of edge defects by embedded eigenvalues of quantum walks}

\author[H. Morioka]{Hisashi MORIOKA}
\address[H. Morioka]{Graduate School of Science and Engineering,
Ehime University, Bunkyo-cho 3, Matsuyama, Ehime, 790-8577, Japan}
\email{morioka@cs.ehime-u.ac.jp}

\author[E. Segawa]{Etsuo SEGAWA}
\address[E. Segawa]{Graduate School Educational Promotion Center, Yokohama National University, Tokiwadai 79-7, Hodogaya-ku, Yokohama 240-8501, Japan}
\email{segawa-etsuo-tb@ynu.ac.jp}

\date{\today}

\begin{document}
\baselineskip 14pt
\maketitle

\begin{abstract}
We consider a position-dependent quantum walk on ${\bf Z}$.
In particular, we derive a detection method for edge defects by embedded eigenvalues of its time evolution operator. 
In the present paper, an edge defect is a set $ \{ y-1 ,y\} $ for $y\in {\bf Z}$ on which the coin operator is an anti-diagonal matrix.  
In fact, under some suitable assumptions, the existence of a finite number of edge defects is equivalent to the existence of embedded eigenvalues of the time evolution operator.
In view of applications, by checking the spectrum, we can detect the existence of disconnecting edge (in the sense of edge defects above) on the line without directly watching it. 

\end{abstract}

%%%%%%%%%%%%%%%%% Section 1 %%%%%%%%%%%%%%%%%%%%

\section{Introduction}

Quantum walks have been studied in various kinds of research fields (see \cite{Am}, \cite{Sh}, \cite{Ven} et al. and its references). 
Recently, there is an abundance of studies on position-dependent quantum walks in view of the spectral theory of unitary operators.
Some results of the weak limit theorem for position-dependent quantum walks were proved by Konno-Luczak-Segawa \cite{KLS}, Endo-Konno \cite{EK1} and Endo et al. \cite{EK2}. 
In view of the scattering theory, the wave operators associated with the time evolution operator were introduced by Suzuki \cite{Su} under the short-range type condition, as well as the asymptotic velocity of the quantum walker and the weak limit theorem were considered as applications. 
We also mention about Richard-Suzuki-Tiedra de Aldecoa \cite{RST}.
A Mourre theory for unitary operators is given and its application to the spectral theory of the quantum walk is derived.

In some models of quantum walks, localization occurs depending on its initial states, and eigenvalues of the time evolution operator have a crucial role in the localization.
If $U$ is a unitary time evolution operator for one-dimensional, two-state quantum walks, eigenvalues and eigenspaces are defined as follows.
If there exists a non-trivial solution $\psi \in \ell^2 ({\bf Z} ; {\bf C}^2 ) $ to the equation $U\psi = e^{i\theta} \psi $ for $ \theta \in [0 , 2\pi )$, we call $ e^{i\theta} $ an eigenvalue of $U$.
Thus the associated eigenspace $\mathcal{E} (\theta )$ is a subspace of $ \ell^2 ({\bf Z} ; {\bf C}^2 )$.
As has been shown by Cantero et al. \cite{CGM}, and Suzuki \cite{Su}, if the initial state has an overlap with $\mathcal{E} (\theta )$ i.e. the initial state is not in $\mathcal{E} (\theta )^{\perp}$ in the sense of $\ell^2 ({\bf Z} ; {\bf C}^2 )$, the localization occurs in the associated quantum walk.
Examples of localizations with one-defect model are in Cantero et al. \cite{CGM}, Konno-Luczak-Segawa \cite{KLS} and Fuda-Funakawa-Suzuki \cite{FuFuSu}.
More generally, we can see a similar result for localizations for quantum walks on graphs (see Segawa-Suzuki \cite{SS}).

In this paper, we consider an approach of \textit{detection of edge defects} by using embedded eigenvalues of the time evolution operator of the one-dimensional, two-state quantum walk.
The rigorous meaning of edge defects will be defined below. 
Let $ \mathcal{H} = \ell^2 ({\bf Z} ; {\bf C}^2 )$ be the space of states. 
The unitary operator $U$ is given by
$$
(U \psi )(x)= P(x+1) \psi (x+1) +Q (x-1) \psi (x-1) , \quad x\in {\bf Z} , 
$$
for every $\psi \in \mathcal{H} $ and 
$$
P(x)= \left[ \begin{array}{cc} 
a (x) & b (x) \\ 0 & 0 \end{array} \right] , \quad Q(x)= \left[ \begin{array}{cc} 
0 & 0 \\ c (x) & d (x) \end{array} \right] .
$$
Here we assume $C(x) := P(x)+Q(x) \in U(2)$ for every $x\in {\bf Z} $ and $U$ is rewritten by $U=SC$ where $S$ is the shift operator defined by
$$
(S\psi )(x)= \left[ \begin{array}{c} \psi^{(0)} (x+1) \\ \psi ^{(1)} (x-1) \end{array} \right] , \quad \psi \in \mathcal{H} , \quad x\in {\bf Z} .
$$
Taking an initial state $\psi_0 \in \mathcal{H}$, we put $\psi (t, \cdot ) := U^t \psi_0$ for $t\in \{ 0,1,2, \ldots \} $.
Since the operator $U$ depends on the position, we call this discrete time evolution \textit{one dimensional position-dependent quantum walk}. 
Thus we call $C$ the \textit{coin operator} of the operator $U$.
The corresponding position-independent quantum walk is given by $U_0 = SC_0 $ where $C_0 := P_0 + Q_0 \in U(2)$ and
$$
P_0 = \left[ \begin{array}{cc} 
a_0  & b_0 \\ 0 & 0 \end{array} \right] , \quad Q_0 = \left[ \begin{array}{cc} 
0 & 0 \\ c_0 & d_0 \end{array} \right] .
$$
We adopt the representation of $C_0$ which is introduced in \cite{RST}.
Precisely, we put $ a_0 = pe^{i\alpha} $, $b_0 = q e^{i\beta} $, $c_0 = -q e^{- i ( \beta - \gamma  )}  $ and $d_0 = p e^{-i ( \alpha - \gamma  )}  $ for $ \alpha  , \beta , \gamma \in [0,2\pi ) $ and $p,q\in [0,1]$ with $p^2 + q^2 =1$ : 
\begin{equation}
C_0 = e^{i\gamma /2} \left[ \begin{array}{cc}
pe^{i(\alpha - \gamma /2 )} & q e^{i (\beta - \gamma /2)} \\ -q e^{-i (\beta - \gamma /2 )} & p e^{-i (\alpha - \gamma /2)} \end{array} \right] .
\label{S1_def_C0}
\end{equation}
Throughout of the paper, we assume that there exist constants $\rho , M>0$ such that 
\begin{equation}
\| C(x) - C_0 \| _{\infty } \leq M e^{-\rho \langle x\rangle } , \quad x\in {\bf Z} ,
\label{S1_eq_exp}
\end{equation}
where $\| \cdot \|_{\infty } $ is the norm of $2\times 2$-matrices defined by 
$$
\| A \| _{\infty} = \max _{1 \leq j,k \leq 2 } | a_{jk} | , \quad A= [a_{jk} ]_{1\leq j,k \leq 2} , 
$$
and $ \langle x \rangle = \sqrt{1+x^2 }$.

%One of attentions of researchers is quantum speed-up of search algorithm (see \cite{Am} and its references.).
In the present paper, we consider the existence or the non-existence of \textit{edge defects} on ${\bf Z}$.
Here we define edge defects as follows.

\begin{definition}
We call the set ${\bf e}_y = \{ y-1 ,y \} $ for $y\in {\bf Z}$ an edge defect if $C(x)=C_1$ for $x\in {\bf e}_y$ where
\begin{equation}
C_1 = e^{ i \gamma ' /2} \left[ \begin{array}{cc}
0 & e^{i (\beta ' - \gamma ' /2)} \\ -e^{-i ( \beta ' - \gamma ' /2)} & 0 \end{array} \right] , 
\label{S1_def_C1}
\end{equation}
for $ \beta ' , \gamma ' \in [0,2\pi )$. 

%We call the set $ {\bf v} _y =\{ y\}$ for $y\in {\bf Z} $ an vertex defect if $C(y)=C_1$ and $a(x) \not= 0$ for $x=y \pm 1 $.
\label{S1_def_edefect}
\end{definition}

Let us make a remark on Definition \ref{S1_def_edefect} in view of applications.
If the edge defect occurs, then there is a disconnection between $ \{y-1,y \} $ in the network
by the definition.
So in this paper we propose a detection way of the existence of a disconnecting part without directly watching it.
Turning our mind to quantum search algorithms driven by quantum walks, we notice that the quantum coins at the target vertices are also perfect reflection operators. 
Then it is possible to interpret that the setting of the edge defect is an {\it infinite system's} analogue of  quantum search algorithms whose target vertices are e.g., $ \{0,1 \} $ ;  in this ``algorithm", we can find how the defects occurs at the targets checking the spectrum of this system (see Figs. \ref{fig_distvertex}-\ref{fig_edge} in \S 5).

Under the assumption (\ref{S1_eq_exp}), we show that one can detect the existence of edge defects by that of eigenvalues of $U$ embedded in the interior of the continuous spectrum $\sigma _{ess} (U)$.
The first result of the present paper is as follows.

\begin{theorem}
Let $p\in (0,1]$.
We assume that there is no edge defect i.e. there exists a constant $\delta >0$ such that $|a(x)| \geq \delta $ for all $x\in {\bf Z} $. 
Moreover, suppose that $C$ and $C_0$ satisfy the condition (\ref{S1_eq_exp}).
Then the continuous spectrum of $U$ is $\sigma_{ess} (U) = \{ e^{i\theta } \ ; \ \theta \in J_{\gamma} \} $ where $J_{\gamma} = J_{\gamma ,1} \cup J_{\gamma ,2} $ with
\begin{gather*}
\begin{split}
& J_{\gamma ,1} = [ \arccos p + \gamma /2 , \pi - \arccos p + \gamma /2 ] , \\
& J_{\gamma ,2} = [ \pi + \arccos  p + \gamma /2 , 2\pi - \arccos p + \gamma /2 ] .
\end{split}
\end{gather*}
Moreover, there is no eigenvalue in $\sigma_{ess} (U) \setminus \mathcal{T} $ where $\mathcal{T} = \{ e^{i\theta} \in \sigma_{ess} (U) \ ; \ \theta \in J_{\gamma , \mathcal{T}} \} $ with
\begin{gather*}
J_{\gamma , \mathcal{T}} =  \left\{ 
\begin{split}
& \arccos p + \gamma /2 ,  \ \pi - \arccos p + \gamma /2 , \\
& \pi + \arccos p + \gamma /2 , \  2\pi - \arccos p + \gamma /2 
\end{split}
\right\} . 
\end{gather*} 
\label{S1_mainthm1}
\end{theorem}

If there are some edge defects, the operator $U$ is given as follows.
Let $C_1 $ be defined by (\ref{S1_def_C1}).
For a positive integer $N>0$, we take $y_1 , \cdots , y_N \in {\bf Z} $, and put
$$
{\bf e} = \bigcup _{j=1}^N {\bf e}_{y_j} ,   \quad {\bf e} _{y_j} = \{ y_j -1 , y_j \} .
$$
For any subset $A\subset {\bf Z} $, let the operator $F_A$ on $\mathcal{H}$ be defined by $ (F_A \psi )(x)= \psi (x)$ for $x\in A$ and $( F_A \psi )(x)=0 $ for $x \in {\bf Z} \setminus A$. 
Then we put 
\begin{equation}
C  =  \sum _{j=1}^N F_{{\bf e}_{y_j}} C_1 + (1-F_{{\bf e}} ) C_2 = F_{{\bf e}} C_1 + (1-F_{{\bf e}} ) C_2 ,
\label{S1_eq_CC}
\end{equation}
where the coin operator $C_2$ given by
$$
 C_2 (x) = \left[ \begin{array}{cc}
a_2 (x) & b_2 (x) \\ c_2 (x) & d_2 (x) \end{array} \right] \in U (2) , \quad x\in {\bf Z} ,
$$
satisfies the assumption (\ref{S1_eq_exp}) and there exists a constant $\delta >0$ such that $|a_2 (x) | \geq \delta $ for all $x\in {\bf Z}$.
In this case, the situation of $U$ and $U_0$ is same as Theorem \ref{S1_mainthm2} in ${\bf Z}\setminus {\bf e}$. 
However, there exists an embedded eigenvalue as follows.

\begin{theorem}
Let $p\in (0,1]$ and $C$ be given by (\ref{S1_eq_CC}). \\
(1) The continuous spectrum of $U$ is $\sigma_{ess} (U)=\{ e^{i\theta} \ ; \ \theta \in J_{\gamma} \} $. \\
(2) For any $\gamma ' \in [0,2\pi )$, we have $\pm i e^{i\gamma ' /2} \in \sigma_p (U)$, and we can take associated eigenfunctions $\Psi _{\pm} \in \mathcal{H}$ such that $\mathrm{supp} \Psi _{\pm} \subset {\bf e}$.  \\
(3) If $ ( \gamma ' + \pi )/2 \in J_{\gamma} \setminus J_{\gamma , \mathcal{T} }$, we have $ \pm i e^{i \gamma ' /2} \in \sigma_p (U) \cap ( \sigma_{ess} (U) \setminus \mathcal{T} )$.
Any associated eigenfunctions $\Psi_{\pm}$ vanish in $\{x\in {\bf Z} \ ; \  x>x^* \, \text{or} \ x<x_* \} $ where $x^* = \max \{ x\in {\bf e} \} $ and $x_* = \min \{ x\in {\bf e} \} $. 
\label{S1_mainthm2}
\end{theorem}

As a consequence of Theorems \ref{S1_mainthm1} and \ref{S1_mainthm2}, we can state the conclusion of this paper.

\begin{cor} 
Let $p\in (0,1]$ and $( \gamma ' + \pi )/2 \in J_{\gamma} \setminus J_{\gamma , \mathcal{T}}$.
Suppose $C$ is given by (\ref{S1_eq_CC}). 
There is no edge defect i.e. $ {\bf e} = \emptyset $ if and only if $U$ has no eigenvalue in $\sigma_{ess} (U)\setminus \mathcal{T}$.
\label{S1_cor_main}
\end{cor}

Theorems \ref{S1_mainthm1} and \ref{S1_mainthm2} are analogues of the Rellich type uniqueness theorem for the Helmholtz equation $(-\Delta -\lambda )u=0$ on the Euclidean space.
Originally it was introduced by Rellich \cite{Re} and Vekoua \cite{Vek}.
This theorem has been generalized to a broad class of partial differential equations, since it plays important roles in the spectral theory (\cite{Tr}, \cite{Li1}, \cite{Li2}, \cite{Ho}, \cite{Mu} and \cite{RaTa}).
Recently, this theorem was generalized for the discrete Schr\"{o}dinger operator on perturbed periodic graphs (\cite{IsMo}, \cite{Ves} and \cite{AIM}).
Note that the Rellich type uniqueness theorem holds in a Banach space larger than $L^2$-space or $\ell^2$-space.
However, it is sufficient to prove in $\ell^2 ({\bf Z} ; {\bf C}^2 )$ for our purpose of the paper.
For the proof, we use a Paley-Wiener theorem and the theory of complex variable.

The plan of this paper is as follows.
In \S 2, we recall basic properties of spectra of unitary operators.
The proof of Theorem \ref{S1_mainthm1} is given in \S 3.
The precise construction of embedded eigenvalues and the associated eigenfunctions are given in \S 4.
We summarize our arguments in \S 5, using some numerical examples.

Throughout of this paper, we use the following basic notations.
We denote the flat torus by ${\bf T} = {\bf R} / ( 2\pi {\bf Z} )$ and the complex torus by ${\bf T} _{{\bf C}} = {\bf C} /( 2\pi {\bf Z} )$.
For any $s\in {\bf R}$, we put $\langle s \rangle = \sqrt{ 1+s^2  }$. 
The unit circle on the complex plane ${\bf C}$ is denoted by $S^1 $.

%%%%%%%%%%%%%%%%%%%%%%%%%%%%%

\section{Continuous spectrum }

\subsection{Spectral decomposition of unitary operators}
Here let us recall some general properties of spectra of unitary operators.
Let $ \mathcal{H} $ be a Hilbert space.
We denote by $ (\cdot , \cdot ) _{\mathcal{H}} $ the inner product of $ \mathcal{H} $ and by $\| \cdot \| _{\mathcal{H}} $ the associated norm.

Let $U$ be a unitary operator on $ \mathcal{H} $.
It is well-known that there exists a spectral decomposition $E_U ( \theta )$ for $ \theta \in {\bf R}$ such that 
$$
U=  \int_0^{2\pi} e^{i\theta} dE_U ( \theta ) ,
$$
where $E_U ( \theta )$ is extended to be zero for $\theta \in (-\infty , 0)$ and to be $1$ for $\theta \in [ 2\pi ,\infty )$.
It is well-known that $\sigma (U) \subset S^1 $.
Since $E_U (\theta )$ is a measure on ${\bf R}$, applying Radon-Nikod\'{y}m theorem, it provides the orthogonal decomposition of $ \mathcal{H} $ associated with $U$ as 
$$
\mathcal{H} = \mathcal{H}_p (U) \oplus \mathcal{H}_{sc} (U) \oplus \mathcal{H} _{ac} (U) ,
$$
where $ \mathcal{H}_p (U ) $ is spanned by eigenfunctions of $U$, $ \mathcal{H}_{sc} (U) $ and $\mathcal{H}_{ac} (U)$ are orthogonal projections on the pure point, the singular continuous and the absolutely continuous subspace of $ \mathcal{H}$, respectively.
Then we put
$$
\sigma_p (U) = \text{the set of eigenvalues of } U \text{ in } \mathcal{H} , 
$$
$$ 
\sigma_{sc} (U)= \sigma ( U | _{\mathcal{H} _{sc} (U)} ) , \quad \sigma_{ac} (U)= \sigma (U | _{\mathcal{H}_{ac} (U)} ) ,
$$
and we call them the point spectrum, the singular continuous spectrum and the absolutely continuous spectrum of $U$, respectively.

We also define the discrete spectrum and the essential spectrum of $U$.
The discrete spectrum $\sigma_d (U)$ is the set of isolated eigenvalues of $U$ with finite multiplicities.
The essential spectrum $\sigma_{ess} (U)$ is defined by $ \sigma_{ess} (U) = \sigma (U) \setminus \sigma_d (U)$. 
Then if $\lambda \in \sigma_{ess} (U)$, $\lambda $ is either an eigenvalue of infinite multiplicity or an accumulation point of $\sigma (U)$.

As in the case of self-adjoint operators, the essential spectrum of $U$ is characterized by singular sequences as follows.

\begin{lemma}
We have $e^{i \theta } \in \sigma_{ess} (U)$ for $ \theta \in [0,2\pi )$ if and only if there exists a sequence $\{ \psi _n \} _{n=1}^{\infty} $ in $ \mathcal{H} $ such that $\| \psi_n \| _{\mathcal{H} } =1$, $\psi_n \to 0$ weakly in $\mathcal{H}$ and $ \| (U-e^{i\theta} )\psi_n \| _{\mathcal{H}} \to 0$ as $n\to \infty $. 
\label{S2_lem_essspec}
\end{lemma}

Proof.
Suppose $ e^{i\theta} \in \sigma_{ess} (U)$.
When $e^{i\theta} $ is an eigenvalue of infinite multiplicities, we can take an orthonormal basis $\{ \psi_n \}_{n=1}^{\infty}  $ in $ \mathrm{Ker} (U-e^{i\theta} )$.
When $e^{i\theta} $ is an accumulation point of $\sigma (U)$, we can take a sequence $ \{ \theta_n \} _{n=1}^{\infty} $ such that $e^{i\theta_n} \in \sigma (U)$ and $\theta_n \to \theta$.
We take sufficiently small $\epsilon _n >0$ so that $I_n = ( \theta_n - \epsilon_n , \theta_n + \epsilon_n )$ satisfies $I_n \cap I_m = \emptyset $ for $m\not= n$.
By choosing $\psi_n \in \mathrm{Ran}( E_U ( I_n ))$ with $\| \psi_n \| _{\mathcal{H}} =1$, we have an orthonormal basis $\{ \psi_n \} _{n=1}^{\infty} $.
Moreover, we obtain 
$$
\| (U-e^{i\theta} ) \psi_n \|^2 _{\mathcal{H}} = \int_{I_n} |e^{is} - e^{i\theta} |^2 d( E_U (s) \psi_n , \psi_n )_{\mathcal{H}} \leq C\epsilon_n^2 \to 0 .
$$

Suppose that there exists a sequence $\{ \psi_n \} _{n=1}^{\infty}$ such that $\psi_n$ satisfies the condition in the statement of the lemma.
If $e^{i\theta} \not\in \sigma (U)$, there exists a constant $\delta >0$ such that $ E_U (( \theta - \delta , \theta + \delta ))=0$ and $ \| (U-e^{i\theta} ) \psi \| _{\mathcal{H}} \geq \delta $ for any $\psi \in \mathcal{H} $.
This is a contradiction. 
If $e^{i \theta } \in \sigma_d (U)$, there exists a constant $ \epsilon >0$ such that $ E_U (( \theta - \epsilon , \theta + \epsilon ))= E_U ( \{ \theta \} )$ for $e^{i\theta} \not= 1$ or $E_U ((-\epsilon , \epsilon )) + E_U ((2\pi - \epsilon , 2\pi + \epsilon ))=E_U (\{ 0\} ) + E_U (\{ 2\pi \} )$ for $e^{i\theta} =1$. 
In the following, we shall prove the case $e^{i\theta} \not= 1$. 
For $e^{i\theta} =1$, the proof is similar.

We can take an orthonormal basis $\{ \phi_j \} _{j=1}^m $ of $\mathrm{Ker} (U-e^{i\theta} )$ for a positive integer $m$. 
Applying the Gram-Schmidt orthonormalization to $\{ \phi_j \} _{j=1}^m \cup \{ \psi _k \} _{k=1}^{\infty }$, we put the resulting sequence $\{ \phi '_j \} _{j=1}^{\infty } $. 
Note that $\phi'_j = \phi_j $ for $j=1 ,\cdots , m$. 
Hence we have $\| (U-e^{i\theta} ) \phi'_j \| _{\mathcal{H}} \to 0 $ as $j\to \infty $.
On the other hand, we have 
$$
\| (U-e^{i\theta} ) \phi '_j \|^2_{\mathcal{H}} = \int _{|s-\theta | \geq \epsilon } |e ^{i (s-\theta )} -1 |^2 d( E_U (s) \phi'_j , \phi '_j )_{\mathcal{H}} \geq \epsilon ^2 ,
$$
for $j>m$. 
This is a contradiction.
\qed

\medskip

As a consequence, we can see that compact perturbations of $U$ do not change its essential spectrum.

\begin{lemma}
Let $U'$ and $U$ be unitary operators on $\mathcal{H}$. 
If $U'-U$ is compact on $\mathcal{H}$, we have $\sigma_{ess} (U') = \sigma_{ess} (U)$. 

\label{S2_lem_weyltype}
\end{lemma}

Proof.
Let $e^{i\theta} \in \sigma_{ess} (U)$.
In view of Lemma \ref{S2_lem_essspec}, there exists a sequence $\{ \psi_n \} _{n=1}^{\infty} $ in $\mathcal{H}$ such that $\| \psi_n \| _{\mathcal{H}} =1$, $\psi_n \to 0$ weakly in $\mathcal{H}$ and $\| (U-e^{i\theta } )\psi_n \| _{\mathcal{H}} \to 0$ as $n\to \infty $. 
Since $U' -U$ is compact, we have $(U'-U )\psi_n \to 0$ in $\mathcal{H}$.
Then we have
$$
\| (U' -e^{i\theta} )\psi_n \| _{\mathcal{H}} \leq \| (U-e^{i\theta } )\psi_n \| _{\mathcal{H}} + \| (U' -U) \psi_n \| _{\mathcal{H}} \to 0.
$$
Applying Lemma \ref{S2_lem_essspec} to $U'$, we obtain $e^{i\theta} \in \sigma_{ess} (U')$.
This implies $ \sigma_{ess} (U) \subset \sigma_{ess} (U') $. 
We can prove $\sigma_{ess} (U' ) \subset \sigma_{ess} (U) $ by the same way.
\qed

%%%%%%%%%%%%%%%%%%%%%%%%%%%%%
\subsection{Essential spectrum}
We turn to the quantum walk.
In the following, the notations $U$ and $U_0$ are used in order to represent the unitary operators of time evolution for the quantum walk, and $\mathcal{H} = \ell^2 ({\bf Z} ; {\bf C}^2 )$.
Let $ \mathcal{F} : \mathcal{H} \to \widehat{\mathcal{H}} := L^2 ({\bf T} ; {\bf C}^2 )$ be the unitary operator defined by 
$$
(\mathcal{F} \psi )(\xi ) = \left[ \begin{array}{c} 
\widehat{\psi }^{(0)} (\xi ) \\ \widehat{\psi }^{(1)} (\xi )
\end{array} \right] ,\quad \widehat{\psi }^{(j)} (\xi )= \frac{1}{\sqrt{2\pi }} \sum _{x\in {\bf Z}} e^{-ix\xi } \psi ^{(j)} (x) , 
$$
for $\xi \in {\bf T}$, $j=0,1$, and every $ \psi \in \mathcal{H} $. 
Then the adjoint operator $ \mathcal{F}^* : \widehat{\mathcal{H}} \to \mathcal{H} $ is given by
$$
(\mathcal{F}^* \widehat{\phi} )(x ) = \left[ \begin{array}{c} 
\phi ^{(0)} (x ) \\ \phi ^{(1)} (x )
\end{array} \right] ,\quad \phi^{(j)} (x) = \frac{1}{\sqrt{2\pi}} \int _{{\bf T}} e^{ix\xi} \widehat{\phi}^{(j)} (\xi ) d\xi ,
$$
for $x\in {\bf Z}$, $j=0,1$, and every $\widehat{\phi} \in \widehat{\mathcal{H}} $.

Letting 
$$
\widehat{U}_0 = \mathcal{F} U_0 \mathcal{F}^* = \mathcal{F} SC_0 \mathcal{F}^*,
$$
we have that $\widehat{U}_0 $ is the operator of multiplication by the unitary matrix 
\begin{equation}
\widehat{U}_0 (\xi ) = \left[ \begin{array}{cc}
a_0 e^{i\xi} & b_0 e^{i\xi} \\ c_0 e^{-i\xi} & d_0 e^{-i\xi} \end{array} \right] .
\label{S2_eq_U0hat}
\end{equation}
In view of (\ref{S1_def_C0}), we have 
\begin{equation}
\widehat{U}_0 (\xi ) =  e^{i\gamma /2}  \left[ \begin{array}{cc}
p e^{i(\alpha - \gamma /2 )} e^{i\xi} & q e^{i(\beta - \gamma /2)} e^{i\xi} \\ -q e^{-i(\beta -\gamma /2)} e^{-i\xi} & p e^{-i( \alpha - \gamma /2 )} e^{-i\xi} \end{array} \right] .
\label{S2_eq_U0hat2}
\end{equation}
Moreover, we obtain for any $ \lambda \in {\bf C}$
\begin{equation}
\mathrm{det} (\widehat{U}_0 (\xi ) -\lambda ) = \lambda^2 -2 \lambda pe^{i\gamma /2} \cos \left( \xi + \alpha - \frac{\gamma}{2} \right) + e^{i\gamma } .
\label{S2_eq_det}
\end{equation}
In view of (\ref{S2_eq_det}), we can see the following fact. 
For the proof, see Lemma 4.1 in \cite{RST}. 

%%%%%%%%%%%%

\begin{lemma}
(1) If $p=0$, we have $ \sigma (U_0 ) = \sigma_p (U_0 ) = \{ \pm i e^{i\gamma /2} \} $. \\
(2) If $p\in (0,1)$, we have $\sigma (U_0 )= \sigma_{ac} (U_0 ) = \{ e^{i \theta} \ ; \ \theta \in J_{\gamma} \}$. \\
(3) If $p=1$, we have $\sigma (U_0 )= \sigma_{ac} (U_0 )= S^1 $.
\label{S2_lem_specU0}
\end{lemma}

In view of the assumption (\ref{S1_eq_exp}), the operator $U-U_0 $ is compact on $\mathcal{H} $.
Applying Lemma \ref{S2_lem_weyltype}, we obtain the following lemma.

\begin{lemma}
(1) If $p\in (0,1) $, we have $\sigma_{ess} (U)=\sigma_{ess} (U_0 ) = \{ e^{i\theta} \ ; \ \theta \in J_{\gamma} \}$. \\
(2) If $p=1$, we have $\sigma_{ess} (U)=\sigma_{ess} (U_0 )  = S^1 $.

\label{S2_lem_essspecUU0}
\end{lemma}

%%%%%%%%%%%%%%%%%%%%%%%%%%%%%%%%%%%%%

\section{Absence of embedded eigenvalues}

\subsection{Thresholds}

Let 
\begin{gather}
M  ( \theta ) = \{ \xi \in {\bf T} \ ; \ p(\xi , \theta )=0  \} ,    \label{S3_def_fermi} \\
M _{reg} ( \theta ) = \{ \xi \in {\bf T} \ ; \ p(\xi , \theta )=0 , \partial _{\xi} p(\xi ,\theta ) \not= 0 \} ,    \label{S3_def_fermi_reg} \\
M_{sng} (\theta ) = \{ \xi \in {\bf T} \ ; \ p(\xi , \theta )=0 , \partial _{\xi} p (\xi , \theta )=0 \} , \label{S3_def_fermising}
\end{gather}
where $p(\xi ,\theta )= \mathrm{det} ( \widehat{U}_0 (\xi ) - e^{i\theta } )$.
Note that $p(\xi ,\theta )$ is a trigonometric polynomial in $\xi $ (see (\ref{S2_eq_det})).
%We put $ \lambda = \kappa_1 + i\kappa_2 $ for $\kappa_1 , \kappa_2 \in {\bf R}$.
%Then we have 
%$$
%p(\xi , \theta )= p_1 (\xi , \theta ) + ip_2 (\xi , \theta ),
%$$
%where 
%\begin{gather*}
%p_1 (\xi , \theta )= -2p \cos \frac{\gamma}{2} \cos \left( \xi + \alpha - \frac{\gamma}{2} \right) +  \kappa_1 ^2 - \kappa_2 ^2 + \cos \gamma , \\
%p_2 (\xi , \theta )= -2p \sin \frac{\gamma}{2} \cos \left( \xi + \alpha - \frac{\gamma}{2} \right) + 2 \kappa_1 \kappa_2 + \sin \gamma .
%\end{gather*}

\begin{lemma}
Suppose $p\in (0,1]$.
If $ \theta \in J_{\gamma} \setminus J_{\gamma ,\mathcal{T}} $, we have $M(\theta )= M_{reg} (\theta )$ and $M_{sng} (\theta )= \emptyset$. 
If $ \theta \in  J_{\gamma ,\mathcal{T}} $, we have $M(\theta )= M_{sng} (\theta )$ and $M_{reg} (\theta )= \emptyset$.
\label{S3_lem_singular}
\end{lemma}

Proof.
Note that
$$
\partial_{\xi} p (\xi , \theta )= 2p e^{i \gamma /2} e^{i\theta} \sin \left( \xi + \alpha  -\frac{\gamma}{2} \right) .
$$
Then $ \partial_{\xi} p (\xi , \theta )=0$ if and only if $ \xi + \alpha - \gamma /2 =0$ modulo $\pi $. 
If $p(\xi , \theta )=\partial_{\xi} p (\xi , \theta ) =0$, we have that $e^{i\theta} $ must be equal to one of the following values : 
$$
 e^{i \gamma /2} \left( p \pm i \sqrt{ 1-p^2 } \right) , \quad  e^{i \gamma /2} \left( - p \pm i \sqrt{ 1-p^2 } \right) .
$$
The lemma follows from these observations.
\qed

%%%%%%%%%%%%%%%%%%%%%%%%%%%%%%%%%%%%%%%
\subsection{Absence of embedded eigenvalues} \label{section_proofmain1}
In \S \ref{section_proofmain1}, we prove Theorem \ref{S1_mainthm1}.
For the proof, we suppose that there exists an eigenvalue in $ \sigma_p (U) \cap ( \sigma_{ess} (U) \setminus \mathcal{T} ) $ and we show a contradiction.

Let us recall the assumptions which we adopt in \S \ref{section_proofmain1} :

\begin{enumerate}

\item $p \in (0,1]$ and there exists a constant $\delta >0$ such that $|a(x)| \geq \delta $ for all $x\in {\bf Z} $.
\item There exist constants $\rho ,M>0$ such that $\| C(x)-C_0 \| _{\infty} \leq Me^{-\rho \langle x\rangle }$ for any $x\in {\bf Z} $.

\end{enumerate}

We assume $e^{i\theta} \in \sigma_p (U) \cap (\sigma_{ess} (U) \setminus \mathcal{T} )$ and let $\psi \in \mathcal{H}$ be the associated eigenfunction.
Putting $f=-(U-U_0 )\psi \in \mathcal{H}$, the equation $(U-e^{i\theta } )\psi =0$ is rewritten as 
$$
(U_0 -e^{i\theta} )\psi =f \quad \text{on} \quad {\bf Z} .
$$
In view of the assumption (2), we have $e^{r \langle \cdot \rangle} f \in \mathcal{H}$ for any $r\in (0,\rho )$.
Passing to the Fourier series, we have 
\begin{equation}
(\widehat{U}_0 (\xi )-e^{i\theta} ) \widehat{\psi} = \widehat{f} \quad \text{on} \quad {\bf  T}.
\label{S3_eq_eigen_torus}
\end{equation}
Moreover, we multiply the equation (\ref{S3_eq_eigen_torus}) by the cofactor matrix of $ \widehat{U}_0 (\xi )-e^{i\theta}$.
Note that each component of the cofactor matrix is trigonometric polynomials.
Then the matrix $ \widehat{U}_0 (\xi )-e^{i\theta}$ is diagonalized and it is sufficient to consider the equation of the form
\begin{equation}
p(\xi , \theta ) \widehat{u} = \widehat{g} \quad \text{on} \quad {\bf  T} ,
\label{S3_eq_eigen_torus2}
\end{equation}
where $\widehat{u} , \widehat{g} \in L^2 ({\bf T})$.

Here we need a Paley-Wiener type theorem.
The following one is Theorem 6.1 in \cite{Ves}.

\begin{theorem}
Let $k_0 >0$ be a constant.  
For a function $\phi \in \ell^2 ({\bf Z} )$, $e^{k\langle \cdot \rangle} \phi \in \ell^2 ({\bf Z} )$ for any $k\in (0,k_0 )$ if and only if the function $\widehat{\phi}$ extends to analytic function in $\{ z\in {\bf T}_{{\bf C}} \ ; \ | \mathrm{Im} \, z| < k_0 / (2\pi ) \}$.
\label{S3_thm_PW}
\end{theorem}

As a direct consequence, we have the following fact.
\begin{lemma}

The function $\widehat{g}$ in (\ref{S3_eq_eigen_torus2}) extends to an analytic function in $\{ z\in {\bf T}_{{\bf C}} \ ; \ | \mathrm{Im} \, z| < \rho / (2\pi ) \}$.
\label{S3_lem_PW}
\end{lemma}

Proof.
Since we have $e^{r \langle \cdot \rangle} f \in \mathcal{H}$ for any $r\in (0, \rho )$, we apply Theorem \ref{S3_thm_PW} to $f$ so that $\widehat{f}$ is analytic in $\{ z\in {\bf T}_{{\bf C}} \ ; \ | \mathrm{Im} \, z| < \rho / (2\pi ) \}$.
Each component of the cofactor matrix is trigonometric polynomials.
Then $\widehat{g}$ is also analytic in $\{ z\in {\bf T}_{{\bf C}} \ ; \ | \mathrm{Im} \, z| < \rho / (2\pi ) \}$.
\qed

\medskip

Next we discuss about the regularity of $\widehat{u}$.

\begin{lemma}
Let $\widehat{u} \in L^2 ({\bf T})$ satisfy the equation (\ref{S3_eq_eigen_torus2}). 
Then $\widehat{u} \in C^{\infty} ({\bf T} )$.
In particular, we have $\widehat{g} (\xi (\theta )) =0$ for $\xi (\theta )\in M(\theta )$. \\

\label{S3_lem_smoothness}
\end{lemma}

Proof.
We take $\xi (\theta ) \in M(\theta )$.
Note that $M(\theta )= M_{reg} (\theta )$ from $e^{i\theta} \in \sigma_p (U) \cap ( \sigma_{ess} (U) \setminus \mathcal{T} )$.
Let $\chi \in C^{\infty} ({\bf T} )$ satisfy $\chi (\xi (\theta )) =1$ with small support. 
In view of $\xi (\theta )\in M_{reg} (\theta )$, we have $ \partial _{\xi} p (\xi (\theta ),\theta ) \not= 0$. 
Thus we can make the change of variable
$$
\eta = \cos \left( \xi + \alpha - \frac{\gamma}{2} \right) - \cos \left( \xi (\theta ) + \alpha - \frac{\gamma}{2} \right) ,
$$
in a small neighborhood of $\xi (\theta )$.
Letting $ \widehat{u} _{\chi} = \chi \widehat{u} $ and $\widehat{g} _{\chi} = \chi \widehat{g} $, we rewrite the equation (\ref{S3_eq_eigen_torus2}) as 
\begin{equation}
\eta   \widehat{u} _{\chi} = -\frac{1}{2p} e^{-i (\theta + \gamma /2)} \widehat{g} _{\chi} \quad \text{on} \quad {\bf T} .
\label{S3_eq_smoothness1}
\end{equation}

Now let us define the Fourier transformation by
$$
\widetilde{u _{\chi}} (t) = \frac{1}{\sqrt{2\pi}} \int_{-\infty}^{\infty} e^{-it\eta} \widehat{u} _{\chi} (\eta )  d\eta , \quad t \in {\bf R} .
$$
We define $\widetilde{g _{\chi}} (t)$ by the same way.
Then the equation (\ref{S3_eq_smoothness1}) is reduced to the differential equation
\begin{equation}
 \partial_t   \widetilde{u_{\chi}} = \frac{i}{2p} e^{-i( \theta + \gamma /2)} \widetilde{g_{\chi}} . 
\label{S3_eq_smoothness2}
\end{equation}
Integrating this equation, we have 
$$
\widetilde{u_{\chi} } (t)= \frac{i}{2p} e^{-i(\theta + \gamma /2)}  \int_0^t \widetilde{g_{\chi}} (s) ds + \widetilde{u_{\chi} } (0)  .
$$
In view of Lemma \ref{S3_lem_PW}, $\widehat{g}_{\chi} $ is smooth. 
Hence $\widetilde{g_{\chi}} $ is rapidly decreasing at infinity.
From $\widehat{u}_{\chi} \in L^2 ({\bf T} )$, we have $\widetilde{u_{\chi}} (t)\to 0$ as $|t| \to \infty $.
Then the limit
$$
\lim _{t\to \infty}  \widetilde{u_{\chi}} (t)= \frac{i}{2p} e^{-i (\theta + \gamma /2)}  \int_0^{\infty}  \widetilde{g_{\chi}} (s) \, ds + \widetilde{u_{\chi}} (0)  ,
$$
exists and we obtain
$$
\widetilde{u_{\chi}} (0)= - \frac{i}{2p} e^{-i (\theta +\gamma /2)} \int_0^{\infty}  \widetilde{g_{\chi} } (s) \, ds .
$$
Therefore, $\widetilde{u_{\chi}} $ is represented by the rapidly decreasing function
\begin{equation}
\widetilde{u_{\chi}} (t)=- \frac{i}{2p} e^{-i (  \theta + \gamma /2)} \int_t^{\infty}  \widetilde{g_{\chi}} (s) \, ds , \quad t \geq 0 .
\label{S3_eq_smoothness3}
\end{equation}
Similarly, we have as $t\to - \infty$ 
$$
\lim _{t\to -\infty}  \widetilde{u_{\chi}} (t)=- \frac{i}{2p} e^{ -i ( \theta + \gamma /2)} \int_{-\infty}^0  \widetilde{g_{\chi}} (s) \, ds + \widetilde{u_{\chi}} (0)  ,
$$
and
$$
\widetilde{u_{\chi}} (0) = \frac{i}{2p} e^{-i (\theta + \gamma /2)}  \int_{-\infty}^0 \widetilde{g_{\chi}} (s) \, ds .
$$
Hence we obtain 
\begin{equation}
\widetilde{u_{\chi}} (t)= \frac{i}{2p} e^{-i (  \theta + \gamma /2)} \int^t_{-\infty}  \widetilde{g_{\chi}} (s) \, ds , \quad t \leq 0 .
\label{S3_eq_smoothness4}
\end{equation}
Then $\widetilde{u_{\chi}} (t)$ is rapidly decreasing as $|t| \to \infty $ and this implies that $\widehat{u}_{\chi} \in C^{\infty} ({\bf T} )$.
Obviously, $\widehat{u}$ is smooth outside any small neighborhood of $\xi (\theta )$.
Then we have $\widehat{u} \in C^{\infty} ({\bf T} )$.
It follows from the equation (\ref{S3_eq_eigen_torus2}) that $\widehat{g} $ vanishes at $\xi (\theta )$. 
\qed

\begin{lemma}
The meromorphic function $ \widehat{g} (z) /p(z, \theta )$ is analytic in $\{ z\in {\bf T} _{{\bf C}} \ ; \ | \mathrm{Im} \, z| < \rho / (2\pi ) \} $.

\label{S3_lem_uextention}
\end{lemma}

Proof.
If $ p(z,\theta )=0$ for $e^{i\theta} \in \sigma_{ess} (U)\setminus \mathcal{T} $, we have
$$
\cos \left( z+\alpha - \frac{\gamma}{2} \right) = \frac{1}{p} \cos \left( \theta -\frac{\gamma}{2} \right) .
$$
This implies $\mathrm{Im} \, z =0$ if $p(z,\theta )=0$ for $e^{i \theta } \in \sigma_{ess} (U) \setminus \mathcal{T} $. 
Therefore, in order to show the analyticity of $\widehat{g} (z) /p(z,\theta )$, it is sufficient to consider a neighborhood of $\xi (\theta ) \in M(\theta )$.
We expand $p(z, \theta )$ and $ \widehat{g} (z) $ into Taylor series at $\xi ( \theta ) \in M(\theta )$ :
$$ 
p(z,\theta )= \sum _{n=0}^{\infty} p_n (z-\xi (\theta ))^n , \quad \widehat{g} (z)= \sum _{n=0}^{\infty} g_n ( z-\xi (\theta ))^n ,
$$
for $p_n , g_n \in {\bf C} $.
In view of $M(\theta )= M_{reg} (\theta )$, we have $p_0 =0$ and $p_1 \not= 0$.
Then Lemma \ref{S3_lem_smoothness} implies $g_0 =0 $ and $\widehat{g} (z) / p(z, \theta )$ is analytic in a neighborhood of $\xi (\theta )$.
The Lemma follows from Lemma \ref{S3_lem_PW}.
\qed

\medskip

In the next step, we show that the eigenfunction $\psi$ decays super-exponentially as $|x| \to \infty $.

\begin{lemma}
For any $k>0$, we have $e^{k \langle \cdot \rangle } \psi \in \mathcal{H} $.
\label{S3_lem_superexp}
\end{lemma}

Proof.
It follows from Lemma \ref{S3_lem_uextention} that the function
$$
u(x ) := \frac{1}{\sqrt{2\pi}} \int _{{\bf T}} e^{ix\xi} \widehat{u} (\xi ) \, d\xi ,
$$
satisfies $e^{r\langle \cdot \rangle} u \in \ell^2 ({\bf Z} )$ for $r\in (0,\rho )$ so that $e^{r \langle \cdot \rangle} \psi \in \mathcal{H} $. 
The assumption (2) implies that the function $f=(U-U_0 )\psi $ satisfies $e^{2r \langle \cdot \rangle } f \in \mathcal{H} $ for any $r\in (0,\rho )$. 
Repeating the arguments in the proofs of Lemmas \ref{S3_lem_PW}-\ref{S3_lem_uextention}, we can see $ e^{2r \langle \cdot \rangle } \psi \in \mathcal{H} $.
We can repeat this procedure any number of times.
Therefore, we have $e^{mr \langle \cdot \rangle} \psi \in \mathcal{H}$ for any $m>0$.
\qed

\medskip

\textit{Proof of Theorem \ref{S1_mainthm1}.}
Plugging Lemmas \ref{S3_lem_PW}-\ref{S3_lem_superexp}, the eigenfunction $\psi $ satisfies $e^{k\langle \cdot \rangle} \psi \in \mathcal{H}$ for any $k>0$. 
The equation $(U-e^{i\theta} )\psi =0$ is rewritten as  
\begin{gather}
a(x+1) \psi ^{(0)} (x+1) + b(x+1) \psi^{(1)} (x+1)=e^{i\theta} \psi ^{(0)} (x) , \label{S3_eq_eigen11} \\
c(x-1) \psi ^{(0)} (x-1) + d(x-1) \psi ^{(1)} (x-1) = e^{i\theta} \psi ^{(1)} (x) . \label{S3_eq_eigen12} 
\end{gather}
Recalling the assumptions (1) and (2), we put
$$
K_1 = \max \left\{ 1 , \, \sup _{x\in {\bf Z}} \| C(x)\| _{\infty} \right\} , \quad K_2 = \max \left\{ 1, \delta^{-1} \right\} .
$$
From the equations (\ref{S3_eq_eigen11}) and (\ref{S3_eq_eigen12}), we have 
\begin{gather*}
\begin{split}
a(x) \psi ^{(0)} (x) = & \, \left( -e^{-i\theta} b(x) c(x-1) + e^{i\theta} \right) \psi ^{(0)} (x-1) \\
& \, - e^{-i\theta} b(x)d(x-1) \psi^{(1)} (x-1) ,
\end{split}
\end{gather*}
and then
$$
| \psi ^{(0)} (x) | \leq 2 K_1^2 K_2 \left( | \psi ^{(0)} (x-1) | + | \psi ^{(1)} (x-1) | \right) .
$$
Repeating the same estimate on the right-hand side, we can see for any $y >0$ that
$$
| \psi ^{(0)} (x) | \leq 2^{2y-1} K_1^{2y} K_2^y \left( | \psi ^{(0)} (x-y) | + | \psi ^{(1)} (x-y) | \right) .
$$
In view of Lemma \ref{S3_lem_superexp}, we obtain 
$$
| \psi ^{(0)} (x) | \leq 2^{2y} K_1^{2y} K_2^y e^{-k\langle x-y \rangle } ,
$$
for any $k>0$. 
Taking a sufficiently large $k$ and tending $y\to \infty$, we see $|\psi^{ (0)} (x)|=0$.
Since $x\in {\bf Z} $ is arbitrary, $\psi^{(0)} $ vanishes on ${\bf Z} $.

Let us go back the equation (\ref{S3_eq_eigen12}).
The equation is rewritten as 
$$
d(x-1) \psi^{(1)} (x-1) = e^{i\theta} \psi ^{(0)} (x) ,
$$
so that
$$
|\psi^{(1)} (x)| \leq K_1 |\psi^{(1)} (x-1)| \leq \cdots \leq K_1^y |\psi^{(1)} (x-y)| ,
$$
for any $y>0$.
Hence we also have
$$
|\psi^{(1)} (x)|\leq K_1^y e^{-k \langle x-y \rangle} ,
$$
for any $k>0$.
Taking a sufficiently large $k>0$ and tending $y\to \infty$, we obtain $\psi ^{(1)} (x)=0$ for any $x\in {\bf Z} $.
\qed

%%%%%%%%%%%%%%%%%%%%%%%%%%%%%%%%

\section{Existence of embedded eigenvalues}
\subsection{Finite support of eigenfunctions}

In this section, we turn to the coin operator $C$ given by (\ref{S1_eq_CC}).
Since $C(x) -C_0$ satisfies the assumption (\ref{S1_eq_exp}), Lemma \ref{S2_lem_essspecUU0} also holds for this case i.e. $\sigma_{ess} (U)= \sigma_{ac} (U_0 )$. 
The set of thresholds $\mathcal{T}$ is also defined by the same manner of Theorem \ref{S1_mainthm1}.
Thus the assertion (1) of Theorem \ref{S1_mainthm2} holds.
On the other hand, the assertion of Theorem \ref{S1_mainthm1} does not hold for this case.
However, we can prove the assertion (3) of Therem \ref{S1_mainthm2} which is weaker than Theorem \ref{S1_mainthm1}.
%This is an analogue of the Rellich type uniqueness theorem (\cite{Re}, \cite{IsMo}, \cite{Ves}, \cite{AIM}).

%\begin{theorem}
%Let $e^{i\theta} \in \sigma_p (U) \cap ( \sigma_{ess} (U) \setminus \mathcal{T} )$ and $\psi \in \mathcal{H} $ be an associated eigenfunction.  
%Then we have $\psi (x)=0$ for $x>x^*$ or $x<x_*$.

%\label{S4_thm_finitecase}
%\end{theorem}

\medskip

\textit{Proof of (3) of Theorem \ref{S1_mainthm2}.} 
We can apply Lemmas \ref{S3_lem_PW}-\ref{S3_lem_superexp} to $U$.
Then we have $e^{k\langle \cdot \rangle} \psi \in \mathcal{H} $ for any $k>0$.
Since we have $a(x)= p e^{i\alpha} \not= 0 $ for $x<x_*$, we can use the estimate which is derived in the proof of Theorem \ref{S1_mainthm1}.
Then we have $\psi =0$ for $x<x_*$. 
In view of the equations (\ref{S3_eq_eigen11}) and (\ref{S3_eq_eigen12}), we have 
\begin{gather*}
\begin{split}
d(x) \psi ^{(1)} (x) = & \, -e^{i\theta} a(x+1) c(x) \psi^{(0)} (x+1)  \\
& \, + \left(  e^{i\theta} - e^{-i\theta} b(x+1)c(x) \right) \psi ^{(1)} (x+1) .
\end{split}
\end{gather*}
Note that $d(x)= pe^{i\alpha} e^{i\gamma} \not= 0$ for $x>x^*$.
Then we have 
$$
| \psi^{(1)} (x) | \leq 2^{2y-1} K_1^{2y} K_2^y e^{-k \langle x+y \rangle} ,
$$
for any large $k>0$ and $y>0$.
We obtain $\psi ^{(0)} (x)=0$ for $x>x^*$ tending $y\to \infty $.
From the equation (\ref{S3_eq_eigen11}), we have
$$
|\psi ^{(0)} (x)|\leq K_1^y | \psi ^{(0)} (x+y)| \leq K_1^y e^{-k \langle x+y\rangle} , 
$$
for any large $k>0$ and $y>0$. 
Hence we also obtain $\psi^{(1)} (x)=0$ for $x>x^*$ tending $y \to \infty $. 
\qed

%%%%%%%%%%%%%%%%%%%%%%%%%%%

\subsection{Embedded eigenvalues}

In order to construct eigenfunctions precisely, we consider the auxiliary operator $U_1 = SC_1$.
Note that $\sigma (U_1 )= \sigma_p (U_1 )= \{ \pm i e^{i\gamma '/2} \} $ (see Lemma \ref{S2_lem_specU0}).

\begin{lemma}
Let $\delta (x) = \delta _{x0}$ for $x\in {\bf Z}$.
Then the function
\begin{equation}
\psi _{\pm} (x)= \frac{1}{\sqrt{2} } \left[ \begin{array}{c}  \mp i e^{i( \beta ' -\gamma '/2 )} \delta (x+1) \\ \delta (x) \end{array} \right] , \quad \beta ' , \gamma ' \in [0,2\pi ),
\label{S4_eq_eigenfunction}
\end{equation}
are normalized eigenfunctions of $U_1$ with eigenvalues $\pm ie^{i\gamma '/2}$, respectively.
\label{S4_lem_supportef}
\end{lemma}

Proof.
The equation $(U_1 - (\pm ie^{i\gamma '/2} )) \psi _{\pm} =0$ is equivalent to
$$
\left[ \begin{array}{cc} 
\mp i e^{i\gamma '/2} & e^{i \beta '} e^{i\xi} \\ -e^{-i\beta '} e^{i\gamma '} e^{-i\xi} & \mp i e^{i\gamma '/2} \end{array} \right]  \left[ \begin{array}{c} \widehat{\psi} _{\pm}^{(0)} (\xi ) \\ \widehat{ \psi }_{\pm}^{(1)} (\xi ) \end{array} \right] =0 , \quad \xi \in {\bf T} .
$$
By a direct computation, we have
$$
\left[ \begin{array}{c} \widehat{\psi} _{\pm}^{(0)} (\xi ) \\ \widehat{ \psi }_{\pm}^{(1)} (\xi ) \end{array} \right]  = s(\xi ) \left[ \begin{array}{c}
\mp i e^{ i(\beta ' -\gamma '/2)} e^{i\xi} \\ 1 \end{array} \right] ,
$$
for any scalar functions $s(\xi )$.
Taking $s(\xi )= (2 \sqrt{\pi} )^{-1}$, we obtain the lemma.
\qed

\medskip

The operator of translation $T_y$ for $y\in {\bf Z}$ is defined by
\begin{equation}
(T_y \psi )(x)= \psi (x-y) , \quad x\in {\bf Z} ,
\label{S4_def_translation}
\end{equation}
for $\psi \in \mathcal{H}$.
Obviously, $T_y \psi _{\pm} $ are also eigenfunctions of $U_1$ with eigenvalues $\pm i e^{i\gamma '/2} $, respectively. 
Moreover, we have $\mathrm{supp} T_y \psi^{(0)} _{\pm} = \{ y-1 \}$ and $ \mathrm{supp} T_y \psi^{(1)} _{\pm} = \{ y \} $.

\medskip

\textit{Proof of (2) of Theorem \ref{S1_mainthm2}.}
We put
$$
\Psi _{\pm} = \kappa_1 T_{y_1} \psi _{\pm}  + \cdots + \kappa_N T_{y_N} \psi _{\pm } ,
$$
for any $ \kappa_1 , \cdots , \kappa_N \in {\bf C} $, where $\psi _{\pm} $ is given by (\ref{S4_eq_eigenfunction}).
Then we have $ \mathrm{supp} \Psi^{(0)}_{\pm} = \{ y_1 -1 , \cdots , y_N -1 \} $ and $\mathrm{supp} \Psi_{\pm}^{(1)} =\{ y_1 , \cdots , y_N \} $. 
Since we have $(F_{{\bf e}_{y_j}} C) \big| _{{\bf e}_{y_j}} = C_1  $ for each $j=1,\cdots ,N$, $\Psi_{\pm} $ satisfies the equation $U\Psi_{\pm} = \pm i e^{i\gamma '/2} \Psi _{\pm} $.
Then $ \pm i e^{i\gamma '/2} \in \sigma_p (U)$ for any $\gamma ' \in [0,2\pi )$.

In view of the assertion (3) of Theorem \ref{S1_mainthm2}, if $ \pm i e^{i\gamma '/2} \in \sigma_p (U) \cap (\sigma_{ess} (U)\setminus \mathcal{T} )$, associated eigenfunctions vanish for $x>x^* $ and $x<x_*$.
\qed

%%%%%%%%%%%%%%%%%%%%%%%%%%%%%

\section{Summary and discussion}

Finally, we summarize our results of the present paper as a conclusive remark by using typical numerical examples.
We consider two typical cases. 
We put $ {\bf e} = {\bf e}_0 \cup {\bf e}_1 = \{ -1 ,0,1 \} $.
Let $U_v = SC_v$ and $U_e = SC_e$ be defined by 
\begin{gather}
C_v =  F_{{\bf e}} \left[ \begin{array}{cc} 1 & 0 \\ 0 & 1 \end{array} \right]   + (1-F_{{\bf e}} ) \left[ \begin{array}{cc} 1/\sqrt{2} & 1/\sqrt{2} \\ -1/\sqrt{2} & 1/\sqrt{2} \end{array} \right]  , \label{S5_def_cv} \\
C_e = F_{{\bf e}}  \left[ \begin{array}{cc} 0 & 1 \\ -1 & 0 \end{array} \right]  + (1-F_{{\bf e}} ) \left[ \begin{array}{cc} 1/\sqrt{2} & 1/\sqrt{2} \\ -1/\sqrt{2} &  1/\sqrt{2} \end{array} \right] . \label{S5_def_ce}
\end{gather} 
For $U_e$, ${\bf e}_0 $ and ${\bf e}_1 $ are edge defects.
On the other hand, $U_v$  does not have edge defects but are perturbed on ${\bf e}$.
From Lemma \ref{S2_lem_essspecUU0}, we have $\sigma_{ess} (U_v)= \sigma_{ess} (U_e) = \{ e^{i\theta} \ ; \ \theta \in J \}$ with  
$$
J= [ \pi /4 , 3\pi /4] \cup [ 5\pi /4 , 7\pi /4 ] .
$$

\begin{figure}[t]
 \begin{minipage}{0.49\hsize}
  \begin{center}
   \includegraphics[width=60mm, bb=0 0 338 229]{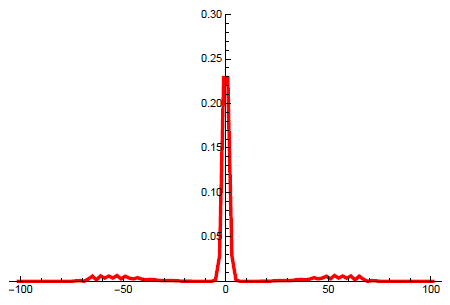}
  \end{center}
  \caption{The distribution of $P_v (X_t =x )$ at $t=100$.}
  \label{fig_distvertex}
 \end{minipage}
 \begin{minipage}{0.49\hsize}
  \begin{center}
   \includegraphics[width=60mm, bb=0 0 336 229]{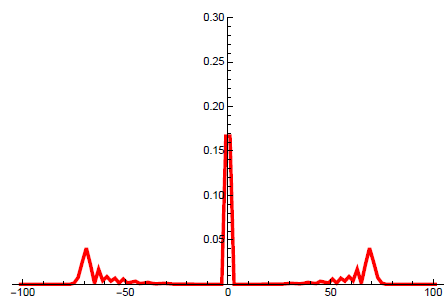}
  \end{center}
  \caption{The distribution of $P_e ( X_t = x)$ at $t=100$.}
  \label{fig_distedge}
\end{minipage}
\end{figure}

Taking the initial state $ \psi_0$ given by 
$$
\psi_0 (x) = \left[ \begin{array}{c} 1/ \sqrt{6} \\ i/\sqrt{6} \end{array} \right] ,  \ x\in {\bf e} , \quad \psi _0 \big| _{{\bf Z} \setminus {\bf e}} =0 ,
$$
we put $\psi_v (t,\cdot ):=U_v^t \psi_0 $ and $\psi_e (t,\cdot ) := U_e^t \psi_0 $ for $t \geq 0$.
Then we compute the probability $P_* (X_t = x)=| \psi_* (t,x)|^2 $ where $*=v$ or $e$ and $X_t$ is the position of the quantum walker at time $t$.
For the numerical results at $t=100$, see Figures \ref{fig_distvertex} and \ref{fig_distedge}.
Localization occurs near $x=0$ for both of $P_v (X_t =x)$ and $P_e (X_t =x)$.
Here localization means $ \limsup _{t\to \infty} P_* (X_t =x) >0$ for some $x\in {\bf Z} $.
Thus we cannot detect edge defects by the existence of localization.

\begin{figure}[t]
 \begin{minipage}{0.49\hsize}
  \begin{center}
   \includegraphics[width=55mm, bb=0 0 323 344]{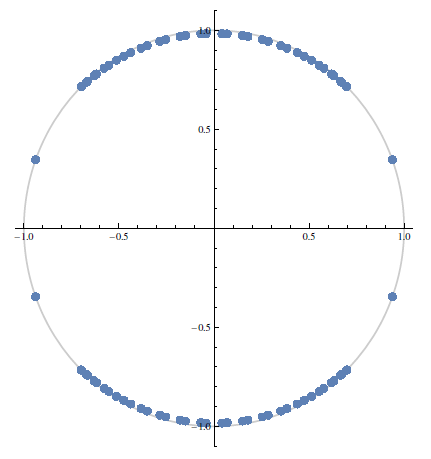}
  \end{center}
  \caption{The distribution of $\sigma (U_v) $.}
  \label{fig_vertex}
 \end{minipage}
 \begin{minipage}{0.49\hsize}
  \begin{center}
   \includegraphics[width=55mm, bb=0 0 323 344]{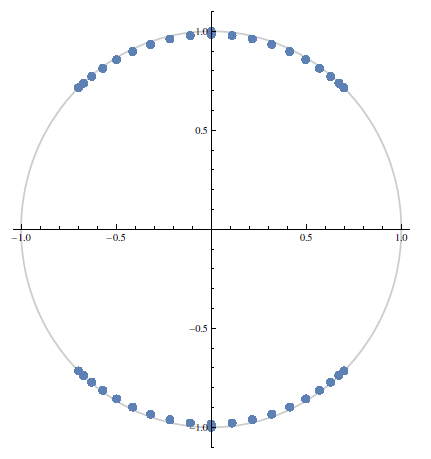}
  \end{center}
  \caption{The distribution of $\sigma (U_e) $.}
  \label{fig_edge}
 \end{minipage}
\end{figure}

If the initial state $ \psi_0 $ has an overlap with an eigenvector of $U_*$, then localization occurs (see \cite{SS}).
For the locations of $\sigma (U_v )$ and $ \sigma ( U_e )$, see Figures \ref{fig_vertex} and \ref{fig_edge}.
$ \sigma _{ess} (U_* )$ is approximated by eigenvalues of the finite rank operator $ U _* \big |_{\{ -60 \leq x \leq 60 \} } $.
The operator $U_v$ has discrete eigenvalues.
On the other hand, $U_e$ has eigenvalues $ \pm i$ which are embedded in the interior of $ \sigma_{ess} (U_e )$.
Localizations of $U_v$ and $U_e$ occur due to eigenvectors of these eigenvalues.
Thus the existence of edge defects is distinguished by the location of eigenvalues.
Precisely, if there exist eigenvalues embedded in the interior of the continuous spectrum, there are some edge defects.

These examples are typical situations to which our main results are applicable (see Theorems \ref{S1_mainthm1} and \ref{S1_mainthm2} and Corollary \ref{S1_cor_main}).

\end{document}